\documentclass[MSNbibl,number,citesort,seceqn,dvips]{arxbj}
\usepackage{upgreek,url,breakurl}
\usepackage{graphicx}
\aid{0}
\volume{20}
\issue{4}
\pubyear{2014}
\firstpage{2305}
\lastpage{2330}
\doi{10.3150/13-BEJ558} %kopijuoti is PTS

\makeatletter

\newcommand{\rright}{\right}
\newcommand{\lleft}{\left}
\newtheorem{theorem}{Theorem}[section]
\newtheorem{corollary}[theorem]{Corollary}
\renewcommand{\pi}{\uppi}
\newcommand{\trup}[2]{{#1}/{#2}}
\newcommand{\C}{\mathbb{C}}
\newcommand{\diag}{\operatorname{diag}}
\newcommand{\tr}{\operatorname{tr}}
\newcommand{\bX}{\mathbf{X}}
\newcommand{\bY}{\mathbf{Y}}
\newcommand{\dd}{\,\mathrm{d}}
\newcommand{\V}{\mathcal{V}}
\newcommand{\E}{\mathbb{E}}
\newcommand{\N}{\mathbb{N}}
\newcommand{\R}{\mathbb{R}}
\makeatother

\begin{document}
\begin{frontmatter}

\title{The affinely invariant distance correlation}
\runtitle{Affinely invariant distance correlation}

\begin{aug}
%%%% inicialai - be tarpu
\author[1]{\inits{J.}\fnms{Johannes} \snm{Dueck}\thanksref{1}},
\author[1]{\inits{D.}\fnms{Dominic} \snm{Edelmann}\thanksref{1}},
\author[2]{\inits{T.}\fnms{Tilmann}~\snm{Gneiting}\thanksref{2}}\\ {\normalfont{and}}
\author[3]{\inits{D.}\fnms{Donald} \snm{Richards}\corref{}\thanksref{3}\ead[label=e4]{richards@stat.psu.edu}}
%%\runauthor{} %% auto

\address[1]{Institut f\"ur Angewandte Mathematik,
Universit\"at Heidelberg, Im Neuenheimer Feld 294, 69120 Heidelberg, Germany}
\address[2]{Heidelberg Institute for Theoretical Studies and Karlsruhe Institute of
Technology, HITS gGmbH, Schloss-Wolfsbrunnenweg 35, 69118 Heidelberg,
Germany}
\address[3]{Department of Statistics,
Pennsylvania State University, University Park, PA 16802, USA.\\
\printead{e4}}
\end{aug}

% HISTORY:
\received{\smonth{4} \syear{2013}}
\revised{\smonth{8} \syear{2013}}

% ABSTRACT
%
\begin{abstract}
Sz\'ekely, Rizzo and Bakirov (\textit{Ann. Statist.} \textbf{35}
(2007) 2769--2794)
and Sz\'ekely and Rizzo (\textit{Ann. Appl. Statist.} \textbf{3}
(2009) 1236--1265), in
two seminal papers, introduced the powerful concept of distance
correlation as a measure of dependence between sets of random
variables. We study in this paper an affinely invariant version of
the distance correlation and an empirical version of that distance
correlation, and we establish the consistency of the empirical quantity.
In the case of subvectors of a multivariate normally distributed random
vector, we provide exact expressions for the affinely invariant distance
correlation in both finite-dimensional and asymptotic settings, and in
the finite-dimensional case we find that the affinely invariant distance
correlation is a function of the canonical correlation coefficients.
To illustrate our results, we consider time series of wind vectors at
the Stateline wind energy center in Oregon and Washington, and we derive
the empirical auto and cross distance correlation functions between wind
vectors at distinct meteorological stations.
\end{abstract}

% KEYWORDS
% visi is mazosios raides ir pagal abecele
%

\end{frontmatter}\vspace*{-3pt}

%s1 #&#
\section{Introduction} \label{sec:introduction}

Sz{\'e}kely, Rizzo and Bakirov \cite{SzeRizBak07} and Sz{\'e}kely
and Rizzo \cite{SzeRiz09},
in two seminal papers, introduced the distance covariance and distance
correlation as powerful measures of dependence. Contrary to the
classical Pearson correlation coefficient, the population distance
covariance vanishes only in the case of independence, and it applies
to random vectors of arbitrary dimensions, rather than to univariate
quantities only.

As noted by Newton \cite{New09}, the ``distance covariance
not only provides a \textit{bona fide} dependence measure, but it does so
with a simplicity to satisfy Don Geman's \textit{elevator test} (i.e., a
method must be sufficiently simple that it can be explained to a colleague
in the time it takes to go between floors on an elevator).'' In
the case of the sample distance covariance, find the pairwise
distances between the sample values for the first variable, and center
the resulting distance matrix; then do the same for the second
variable. The square of the sample distance covariance equals the
average entry in the componentwise or Schur product of the two
centered distance matrices. Given the theoretical appeal of the
population quantity, and the striking simplicity of the sample version,
it is not surprising that the distance covariance is experiencing a
wealth of applications, despite having been introduced merely half a
decade ago.\vadjust{\goodbreak}

Specifically, let $p$ and $q$ be positive integers. For column
vectors $s \in\R^p$ and $t \in\R^q$, denote by $|s|_p$ and $|t|_q$
the standard Euclidean norms on the corresponding spaces; thus, if $s
= (s_1,\ldots,s_p)'$ then
\[
|s|_p = \bigl(s_1^2+\cdots+s_p^2
\bigr)^{1/2},
\]
and similarly for $|t|_q$. For vectors $u$ and $v$ of the same
dimension, $p$, we let $\langle u,v \rangle_p$ be the standard
Euclidean scalar product of $u$ and $v$. For jointly distributed
random vectors $X \in\R^p$ and $Y \in\R^q$, let
\[
f_{X,Y}(s,t) = \E\exp \bigl[ \mathrm{i} \langle s,X
\rangle_p + \mathrm{i} \langle t,Y \rangle_q \bigr]
\]
be the joint characteristic function of $(X,Y)$, and let $f_X(s) =
f_{X,Y}(s,0)$ and $f_Y(t) = f_{X,Y}(0,t)$ be the marginal
characteristic functions of $X$ and $Y$, where $s \in\R^p$ and $t \in
\R^q$. Sz\'ekely \textit{et al.} \cite{SzeRizBak07} introduced the
\emph{distance covariance} between $X$ and $Y$ as the
nonnegative number $\V(X,Y)$ defined by
%
%e1.1 #&#
%
\begin{equation}
\label{eq:dcov} \V^2(X,Y) = \frac{1}{c_p c_q} \int
_{\R^{p+q}} \frac{|f_{X,Y}(s,t)-f_X(s)f_Y(t)|^2} {
|s|_p^{p+1} |t|_q^{q+1}} \dd s \dd t,
\end{equation}
where $|z|$ denotes the modulus of $z \in\C$ and
%
%e1.2 #&#
%
\begin{equation}
\label{eq:cp} c_p = \frac{\pi^{(\trup{1}{2})(p+1)}}{\Gamma ((\trup
{1}{2})(p+1) )}.
\end{equation}
The \emph{distance correlation} between $X$ and $Y$ is the
nonnegative number defined
by
%
%e1.3 #&#
%
\begin{equation}
\label{eq:dcor} \mathcal{R}(X,Y) = \frac{\V(X,Y)}{\sqrt{\V(X,X)\V(Y,Y)}}
\end{equation}
if both $\V(X,X)$ and $\V(Y,Y)$ are strictly positive, and defined
to be zero otherwise. For distributions with finite first moments,
the distance correlation characterizes independence in that $0 \le
\mathcal{R}(X,Y) \le1$ with $\mathcal{R}(X,Y) = 0$ if and only if $X$
and $Y$ are independent.

A crucial property of the distance correlation is that it is invariant
under transformations of the form
%
%e1.4 #&#
%
\begin{equation}
\label{eq:dcor.invariance} (X,Y) \longmapsto(a_1 + b_1
C_1 X,a_2 + b_2 C_2 Y),
\end{equation}
where $a_1 \in\R^p$ and $a_2 \in\R^q$, $b_1$ and $b_2$ are nonzero
real numbers, and the matrices $C_1 \in\R^{p \times p}$ and
$C_2 \in\R^{q \times q}$ are orthogonal. However, the distance
correlation fails to be invariant under the group of all invertible
affine transformations of
$(X,Y)$, which led
Sz\'ekely \textit{et al.} \cite{SzeRizBak07}, pages 2784--2785, and
Sz\'ekely and Rizzo \cite{SzeRiz09}, pages 1252--1253, to propose
an affinely
invariant sample version of the distance correlation.

Adapting this proposal to the population setting, the
\emph{affinely invariant distance covariance} between distributions
$X$ and $Y$ with finite second moments and nonsingular population
covariance matrices ${\Sigma_{X}}$ and ${\Sigma_{Y}}$, respectively,
can be introduced
as the nonnegative number $\widetilde\V(X,Y)$ defined by
%
%e1.5 #&#
%
\begin{equation}
\label{eq:aidcov.def} \widetilde\V^2(X,Y) = \V^2\bigl(
\Sigma_X^{-1/2}X,\Sigma_Y^{-1/2}Y\bigr).
\end{equation}
The \emph{affinely invariant distance correlation} between
$X$ and $Y$ is the nonnegative number defined by
%
%e1.6 #&#
%
\begin{equation}
\label{eq:aidcor.def} \widetilde{\mathcal{R}}(X,Y) =
\frac{\widetilde\V(X,Y)}{\sqrt{\widetilde\V(X,X) \widetilde\V(Y,Y)\vphantom{\sum}}}
\end{equation}
if both $\widetilde\V(X,X)$ and $\widetilde\V(Y,Y)$ are strictly
positive, and defined to be zero otherwise. In the sample versions
proposed by Sz\'ekely \textit{et al.} \cite{SzeRizBak07}, the
population quantities are
replaced by their natural estimators. Clearly, the population
affinely invariant distance correlation and its sample version are
invariant under the group of invertible affine transformations, and in
addition to satisfying this often-desirable group invariance property
(Eaton \cite{Eat89}), they inherit the desirable
properties of the standard
distance dependence measures. In particular,
$0 \le\widetilde{\mathcal{R}}(X,Y) \leq1$ and, for populations with
finite second moments and positive definite covariance matrices,
$\widetilde{\mathcal{R}}(X,Y) = 0$ if and only if $X$ and $Y$ are
independent.

The remainder of the paper is organized as follows. In Section~\ref
{sec:sample}, we review the sample version of the affinely
invariant distance correlation introduced by Sz\'ekely \textit{et al.}
\cite{SzeRizBak07},
and we prove that the sample version is strongly consistent. In
Section~\ref{sec:aidc.mvn}, we provide exact expressions for the
affinely invariant distance correlation in the case of subvectors from
a multivariate normal population of arbitrary dimension, thereby
generalizing a result of Sz\'ekely \textit{et al.} \cite{SzeRizBak07}
in the bivariate
case; our result is non-trivial, being derived using the theory of
zonal polynomials and the hypergeometric functions of matrix argument,
and it enables the explicit and efficient calculation of the affinely
invariant distance correlation in the multivariate normal case.

In Section~\ref{sec:limits}, we study the behavior of the affinely
invariant distance measures for subvectors of multivariate normal
populations in limiting cases as the Frobenius norm of the
cross-covariance matrix converges to zero, or as the
dimensions of the subvectors converge to infinity. We expect that
these results will motivate and provide the theoretical basis for many
applications of distance correlation measures for high-dimensional
data.

As an illustration of our results, Section~\ref{sec:Stateline}
considers time series of wind vectors at the Stateline wind energy
center in Oregon and Washington; we shall derive the empirical auto
and cross distance correlation functions between wind vectors at
distinct meteorological stations. Finally, we provide in Section~\ref
{sec:discussion} a discussion in which we make a case for the
use of the distance correlation and the affinely invariant distance
correlation, which we believe to be appealing and powerful
multivariate measures of dependence.

%s2 #&#
\section{The sample version of the affinely invariant distance correlation}
\label{sec:sample}

In this section, which is written primarily to introduce readers
to distance correlation measures, we
describe sample versions of the affinely invariant distance covariance
and distance correlation as introduced by Sz\'ekely \textit{et al.}
\cite{SzeRizBak07},
pages 2784--2785, and Sz\'ekely and Rizzo \cite{SzeRiz09}, pages 1252--1253.

First, we review the sample versions of the standard distance
covariance and distance correlation. Given a random sample
$(X_1,Y_1), \ldots, (X_n,Y_n)$ from jointly distributed random
vectors $X \in\R^p$ and $Y \in\R^q$, we set
\[
\bX= [ X_1, \ldots, X_n ] \in\R^{p \times n} \quad
\mbox{and}\quad \bY= [ Y_1, \ldots, Y_n ] \in
\R^{q \times n}.
\]
A natural way of introducing a
sample version of the distance covariance is to let
\[
f_{\bX,\bY}^n(s,t) = \frac{1}n \sum
_{j=1}^n \exp \bigl[ \mathrm{i} \langle
s,X_j \rangle_p + \mathrm{i} \langle t,Y_j
\rangle_q \bigr]
\]
be the corresponding empirical characteristic function, and to write
$f^n_\bX(s) = f_{\bX,\bY}^n(s,0)$ and $f^n_\bY(t) = f_{\bX,\bY}^n(0,t)$
for the respective marginal empirical characteristic functions. The
\emph{sample distance covariance} then is the nonnegative number
$\V_n(\bX,\bY)$ defined by
\[
\V_n^2(\bX,\bY) = \frac{1}{c_p c_q} \int
_{\R^{p+q}} \frac{|f_{\bX,\bY}^n(s,t)-f_\bX^n(s)f_\bY^n(t)|^2}{|s|_p^{p+1}
|t|_q^{q+1}} \dd s \dd t,
\]
where $c_p$ is the constant given in (\ref{eq:cp}).

Sz\'ekely \textit{et al.} \cite{SzeRizBak07}, in a \emph{tour de force},
showed that
%
%e2.1 #&#
%
\begin{equation}
\label{eq:sample.dcov} \V_n^2(\bX,\bY) = \frac{1}{n^2} \sum
_{k,l=1}^n A_{kl} B_{kl},
\end{equation}
where
\[
a_{kl} = |X_k-X_l|_p,\qquad
\bar{a}_{k\boldsymbol\cdot} = \frac{1}{n} \sum_{l=1}^n
a_{kl},\qquad \bar{a}_{\boldsymbol\cdot l} = \frac{1}{n} \sum
_{k=1}^n a_{kl},\qquad
\bar{a}_{\boldsymbol\cdot\boldsymbol\cdot} = \frac{1}{n^2} \sum_{k,l=1}^n
a_{kl}
\]
and
\[
A_{kl} = a_{kl} - \bar{a}_{k\boldsymbol\cdot} - \bar
{a}_{\boldsymbol\cdot l} + \bar{a}_{\boldsymbol\cdot\boldsymbol\cdot},
\]
and similarly for $b_{kl} = |Y_k-Y_l|_q$,
$\bar{b}_{k\boldsymbol\cdot}$, $\bar{b}_{\boldsymbol\cdot l}$,
$\bar{b}_{\boldsymbol\cdot\boldsymbol\cdot}$, and $B_{kl}$, where
$k,l = 1,\ldots,n$. Thus, the squared sample distance covariance
equals the average entry in the componentwise or Schur product of the
centered distance matrices\vadjust{\goodbreak} for the two variables. The
\emph{sample distance correlation} then is defined by
%
%e2.2 #&#
%
\begin{equation}
\label{eq:sample.dcor} \mathcal{R}_n(\bX,\bY) = \frac{\V_n(\bX,\bY)}{\sqrt{\V_n(\bX,\bX)\V_n(\bY,\bY)}}
\end{equation}
if both $\V_n(\bX,\bX)$ and $\V_n(\bY,\bY)$ are strictly
positive, and
defined to be zero otherwise. Computer code for calculating these sample
versions is available in an \textsc{R} package by Rizzo and
Sz{\'e}kely~\cite{autokey16}.

Now let $S_\bX$ and $S_\bY$ denote the usual sample covariance matrices
of the data $\bX$ and $\bY$, respectively. Following Sz\'ekely
\textit{et al.} \cite{SzeRizBak07}, page 2785, and Sz\'ekely and Rizzo
\cite{SzeRiz09},
page 1253, the
\emph{sample affinely invariant distance covariance} is the nonnegative
number $\widetilde\V_n(\bX,\bY)$ defined by
%
%e2.3 #&#
%
\begin{equation}
\label{eq:sample.aidcov} \widetilde\V_n^2(\bX,\bY) =
\V_n^2\bigl( S_\bX^{-1/2}\bX,
S_\bY ^{-1/2}\bY\bigr)
\end{equation}
if $S_\bX$ and $S_\bY$ are positive definite, and defined to be zero
otherwise. The \emph{sample affinely invariant distance correlation}
is defined by
%
%e2.4 #&#
%
\begin{equation}
\label{eq:sample.aidcor} \widetilde{\mathcal{R}}_n(\bX,\bY) =
\frac{\widetilde{\V}_n(\bX,\bY)} {
\sqrt{\widetilde{\V}_n(\bX,\bX) \widetilde\V_n(\bY,\bY)\vphantom{\sum}}}
\end{equation}
if the quantities in the denominator are strictly positive, and defined
to be zero otherwise. The sample affinely invariant distance
correlation inherits the properties of the sample distance
correlation; in particular
\[
0 \le\widetilde{\mathcal{R}}_n(\bX,\bY) \le1,
\]
and $\widetilde{\mathcal{R}}_n(\bX,\bY) = 1$ implies that $p = q$, that
the linear spaces spanned by $\bX$ and $\bY$ have full rank, and that
there exist a vector $a \in\R^p$, a nonzero number $b \in\R$, and
an orthogonal
matrix $C \in\R^{p \times p}$ such that $S_\bY^{-1/2} \bY= a + b C
S_\bX^{-1/2} \bX$.

Our next result shows that the sample affinely invariant distance
correlation is a consistent estimator of the respective population quantity.

%th2.1 #&#
%
\begin{theorem} \label{th:consistency}
Let $(X,Y) \in\R^{p+q}$ be jointly distributed random vectors with
positive definite marginal covariance matrices
$\Sigma_X \in\R^{p \times p}$ and $\Sigma_Y \in\R^{q \times q}$,
respectively. Suppose that $(X_1,Y_1),\ldots,(X_n,Y_n)$ is a random
sample from $(X,Y)$, and let $\bX= [X_1,\ldots,X_n] \in\R^{p \times n}$
and $\bY= [Y_1,\ldots,Y_n] \in\R^{q \times n}$. Also,
let $\widehat{\Sigma}_\bX$ and $\widehat{\Sigma}_\bY$ be strongly
consistent estimators for ${\Sigma_{X}}$ and ${\Sigma_{Y}}$,
respectively. Then
\[
\V_n^2 \bigl(\widehat{\Sigma}_\bX^{-1/2}
\bX, \widehat{\Sigma}_\bY ^{-1/2} \bY\bigr) \to\widetilde{
\V}^2(X,Y),
\]
almost surely, as $n \to\infty$. In particular, the sample affinely
invariant distance correlation satisfies
%
%e2.5 #&#
%
\begin{equation}
\label{eq:consistency} \widetilde{\mathcal{R}}_n(\bX,\bY) \to\widetilde{
\mathcal{R}}(X,Y),
\end{equation}
almost surely.
\end{theorem}

\begin{pf}
As the covariance matrices ${\Sigma_{X}}$ and ${\Sigma_{Y}}$ are
positive definite, we
may assume that the strongly consistent estimators
$\widehat{\Sigma}_\bX$ and $\widehat{\Sigma}_\bY$ also are positive
definite. Therefore, in order to prove the first statement it suffices
to show that
%
%e2.6 #&#
%
\begin{equation}
\label{eq:proof1} \V_n^2 \bigl(\widehat{
\Sigma}_\bX^{-1/2} \bX, \widehat{\Sigma}_\bY
^{-1/2} \bY\bigr) - \V_n^2 \bigl(
\Sigma_X^{-1/2} \bX, \Sigma_Y^{-1/2} \bY
\bigr) \to0,
\end{equation}
almost surely. By the decomposition of Sz\'ekely \textit{et al.} \cite
{SzeRizBak07},
page 2776, equation (2.18), the left-hand side of (\ref{eq:proof1})
can be written as an average of terms of the form
\[
\bigl | \widehat{\Sigma}_\bX^{-1/2}(X_k-X_l)
\bigr |_p \bigl | \widehat{\Sigma}_\bY^{-1/2}(Y_k-Y_m)
\bigr |_q - \bigl | \Sigma_X^{-1/2}(X_k-X_l)
\bigr |_p \bigl | \Sigma_Y^{-1/2}(Y_k-Y_m)
\bigr |_q.
\]
Using the identity
\begin{eqnarray*}
&& \bigl | \widehat{\Sigma}_\bX^{-1/2}(X_k-X_l)
\bigr |_p \bigl | \widehat{\Sigma}_\bY^{-1/2}(Y_k-Y_m)
\bigr |_q
\\
&&\quad= \bigl | \bigl(\widehat{\Sigma}_\bX^{-1/2} -
\Sigma_X^{-1/2} + \Sigma_X^{-1/2}\bigr)
(X_k-X_l) \bigr |_p \bigl | \bigl(\widehat{
\Sigma}_\bY^{-1/2} - \Sigma_Y^{-1/2} +
\Sigma_Y^{-1/2}\bigr) (Y_k-Y_m)
\bigr |_q,
\end{eqnarray*}
we obtain
\begin{eqnarray*}
&& \bigl | \widehat{\Sigma}_\bX^{-1/2}(X_k-X_l)
\bigr |_p \bigl | \widehat{\Sigma}_\bY^{-1/2}(Y_k-Y_m)
\bigr |_q - \bigl | \Sigma_X^{-1/2}(X_k-X_l)
\bigr |_p \bigl | \Sigma_Y^{-1/2}(Y_k-Y_m)
\bigr |_q
\\
&&\quad \le\bigl \| \widehat{\Sigma}_\bX^{-1/2} -
\Sigma_X^{-1/2} \bigr \| \bigl \| \widehat{\Sigma}_\bY^{-1/2}
- \Sigma_Y^{-1/2} \bigr \| |X_k-X_l|_p
|Y_k-Y_m|_q
\\
&&\qquad{} + \bigl \| \widehat{\Sigma}_\bX^{-1/2} -
\Sigma_X^{-1/2} \bigr \| | X_k-X_l|_p
\bigl | \Sigma_Y^{-1/2}(Y_k-Y_m)
\bigr |_q
\\
&&\qquad{} + \bigl \| \widehat{\Sigma}_\bY^{-1/2} -
\Sigma_Y^{-1/2} \bigr \|\bigl  | \Sigma_X^{-1/2}(X_k-X_l)
\bigr |_p |Y_k-Y_m|_q,
\end{eqnarray*}
where the matrix norm $\|\Lambda\|$ is the largest eigenvalue of
$\Lambda$ in absolute value. Now we can separate the three sums in
the decomposition of Sz\'ekely \textit{et al.} \cite{SzeRizBak07},
page 2776,
equation
(2.18) and place the factors like $\| \widehat{\Sigma}_\bX^{-1/2} -
\Sigma_{X}^{-1/2} \|$ in front of the sums, since they appear in every
summand. Then, $\| \widehat{\Sigma}_\bX^{-1/2} - \Sigma_X^{-1/2} \|$
and $\| \widehat{\Sigma}_\bY^{-1/2} - \Sigma_Y^{-1/2} \|$ tend to zero
and the remaining averages converge to constants (representing some
distance correlation components) almost surely as $n \to\infty$, and
this completes the proof of the first statement. Finally, the
property (\ref{eq:consistency}) of strong consistency of
$\widetilde{\mathcal{R}}_n(\bX,\bY)$ is obtained immediately upon
setting $\widehat{\Sigma}_\bX= S_\bX$ and
$\widehat{\Sigma}_\bY= S_\bY$.
\end{pf}

Sz\'ekely \textit{et al.} \cite{SzeRizBak07}, page 2783, proposed a
test for independence that
is based on the sample distance correlation. From their results, we
see that the asymptotic properties of the test statistic are not
affected by the transition from the standard distance correlation to
the affinely invariant distance correlation. Hence, a completely
analogous but different test can be stated in terms of the affinely
invariant distance correlation. Noting the results of Kosorok
\cite{Kos09}, Section~4; \cite{Kos13}, we raise the possibility that the
specific details
can be devised in a judicious, data-dependent way so that the power
of the test for independence increases when the transition is made to
the affinely invariant distance correlation. Alternative multivariate tests
for independence based on distances have recently been proposed by Heller
\textit{et al.} \cite{HelHelGor13} and Sz{\'e}kely and Rizzo \cite{SzeRiz13}.

%s3 #&#
\section{The affinely invariant distance correlation for multivariate
normal populations}
\label{sec:aidc.mvn}

We now consider the problem of calculating the affinely invariant
distance correlation between the random vectors $X$ and $Y$ where
$(X,Y) \sim\mathcal{N}_{p+q}(\mu,\Sigma)$,
a multivariate normal distribution with mean vector $\mu\in\R^{p+q}$,
covariance matrix $\Sigma\in\R^{(p+q) \times(p+q)}$, where $X$ and
$Y$ have nonsingular marginal covariance matrices
${\Sigma_{X}}\in\R^{p \times p}$ and ${\Sigma_{Y}}\in\R^{q \times
q}$, respectively.

For the case in which $p = q = 1$, that is, the bivariate normal
distribution, the problem was solved by Sz\'ekely \textit{et al.}
\cite{SzeRizBak07}. In
that case, the formula for the affinely invariant distance correlation
depends only on $\rho$, the correlation coefficient, and appears in
terms of the functions $\sin^{-1}\rho$ and $(1-\rho^2)^{1/2}$, both of
which are well-known to be special cases of Gauss' hypergeometric
series. Therefore, it is natural to expect that the general case will
involve generalizations of Gauss' hypergeometric series, and
Theorem~\ref{th:aidcov} below demonstrates that such is indeed the
case. To formulate this result, we need to recall the rudiments of the
theory of zonal polynomials (Muirhead \cite{Mui82}, Chapter~7).

A \emph{partition} $\kappa$ is a vector of nonnegative integers
$(k_1,\ldots,k_q)$ such that $k_1 \ge\cdots\ge k_q$. The integer
$|\kappa| = k_1 + \cdots+ k_q$ is called the \emph{weight} of
$\kappa$; and $\ell(\kappa)$, the \emph{length} of $\kappa$, is the
largest integer $j$ such that $k_j > 0$. The \emph{zonal polynomial}
$C_\kappa(\Lambda)$ is a polynomial mapping from the class
of symmetric matrices $\Lambda\in\R^{q \times q}$ to the real line
which satisfies several properties, the following of which are crucial
for our results:
\begin{enumerate}[(b)]
\item[(a)]
Let $O(q)$ denote the group of orthogonal matrices in $\R^{q \times q}$.
Then
%
%e3.1 #&#
%
\begin{equation}
\label{eq:property.a} C_\kappa\bigl(K' \Lambda K\bigr) =
C_\kappa(\Lambda)
\end{equation}
for all $K \in O(q)$; thus, $C_\kappa(\Lambda)$ is a symmetric
function of
the eigenvalues of $\Lambda$.

\item[(b)]
The polynomial $C_\kappa(\Lambda)$ is homogeneous of degree $|\kappa|$
in $\Lambda$: For any $\delta\in\R$,
%
%e3.2 #&#
%
\begin{equation}
\label{eq:homogeneous} C_\kappa(\delta\Lambda) = \delta^{|\kappa|}
C_\kappa(\Lambda).
\end{equation}

\item[(c)]
If $\Lambda$ is of rank $r$, then $C_\kappa(\Lambda) = 0$ whenever
$\ell(\kappa) > r$.

\item[(d)]
For any nonnegative integer $k$,
%
%e3.3 #&#
%
\begin{equation}
\label{eq:property.d} \sum_{|\kappa|=k} C_\kappa(\Lambda)
= (\tr\Lambda)^k.
\end{equation}

\item[(e)]
For any symmetric matrices $\Lambda_1, \Lambda_2 \in\R^{q \times q}$,
%
%e3.4 #&#
%
\begin{equation}
\label{eq:property.e} \int_{O(q)} C_\kappa
\bigl(K'\Lambda_1 K \Lambda_2\bigr) \dd K =
\frac{C_\kappa(\Lambda_1) C_\kappa(\Lambda_2)}{C_\kappa(I_q)},
\end{equation}
where $I_q = \diag(1, \ldots, 1) \in\R^{q \times q}$ denotes the
identity matrix and the integral is with respect to the Haar measure
on $O(q)$, normalized to have total volume 1.

\item[(f)] Let $\lambda_1,\ldots,\lambda_q$ be the eigenvalues
of $\Lambda$. Then, for a partition $(k)$ with one part,
%
%e3.5 #&#
%
\begin{equation}
\label{eq:zonal.onepart} C_{(k)}(\Lambda) = \frac{k!}{(\trup{1}{2})_k} \sum
_{i_1+\cdots
+i_q = k} \prod_{j=1}^q
\frac{(\trup{1}{2})_{i_j} \lambda_j^{i_j}}{i_j!},
\end{equation}
where the sum is over all nonnegative integers $i_1,\ldots,i_q$ such
that $i_1+\cdots+i_q = k$, and
\[
(\alpha)_k = \frac{\Gamma(\alpha+k)}{\Gamma(\alpha)} = \alpha(\alpha+1) (\alpha+2)
\cdots(\alpha+k-1),
\]
$\alpha\in\C$, is standard notation for the rising factorial. In
particular, on setting $\lambda_j = 1$, $j=1,\ldots,q$, we obtain
from (\ref{eq:zonal.onepart})
%
%e3.6 #&#
%
\begin{equation}
\label{eq:Ck} C_{(k)}(I_q) = \frac{((\trup{1}{2}) q)_k}{(\trup{1}{2})_k}
\end{equation}
(Muirhead \cite{Mui82}, page 237, equation (18),
Gross and Richards \cite{GroRic87},
page 807, Lemma~6.8).
\end{enumerate}

With these properties of the zonal polynomials, we are ready to
state our key result which obtains an explicit formula for the
affinely invariant distance covariance in the case of a Gaussian
population of arbitrary dimension and arbitrary covariance matrix
with positive definite marginal covariance matrices. This formula
turns out to be a function depending only on the dimensions $p$ and $q$
and the eigenvalues of the matrix
$\Lambda= \Sigma_Y^{-1/2} {\Sigma_{YX}}\Sigma_X^{-1} {\Sigma
_{XY}}\Sigma_Y^{-1/2}$,
that is, the squared canonical correlation coefficients of the
subvectors $X$ and $Y$.
For fixed dimensions this implies
$\widetilde{\mathcal{R}}(X,Y)=g(\lambda_1,\ldots,\lambda_r)$, where
$r=\min(p,q)$ and $\lambda_1,\ldots,\lambda_r$ are the canonical correlation
coefficients of $X$ and $Y$. Due to the functional invariance, the
maximum likelihood estimator (MLE) for the affinely invariant distance
correlation in the Gaussian setting
is hence defined by $g(\widehat{\lambda}_1,\ldots,\widehat{\lambda
}_r)$, where
$\widehat{\lambda}_1,\ldots,\widehat{\lambda}_r$ are the MLEs of the
canonical correlation coefficients.

%th3.1 #&#
%
\begin{theorem} \label{th:aidcov}
Suppose that $(X,Y) \sim\mathcal{N}_{p+q}(\mu,\Sigma)$, where
\[
\Sigma= \lleft( %
\matrix{ {\Sigma_{X}}& {
\Sigma_{XY}}
\cr
{\Sigma_{YX}}& {\Sigma_{Y}} }
\rright)
\]
with ${\Sigma_{X}}\in\R^{p \times p}$, ${\Sigma_{Y}}\in\R^{q
\times q}$, and
${\Sigma_{XY}}\in\R^{p \times q}$. Then
%
%e3.7 #&#
%
\begin{equation}
\label{eq:aidcov.mvn} \widetilde{\V}^2(X,Y) = 4\pi\frac{c_{p-1}}{c_p}
\frac{c_{q-1}}{c_q} \sum_{k=1}^\infty
\frac{2^{2k}-2}{k! 2^{2k}} \frac{(\trup{1}{2})_k (-\trup{1}{2})_k (-\trup{1}{2})_k} {
((\trup{1}{2}) p)_k ((\trup{1}{2}) q)_k} C_{(k)}(\Lambda),
\end{equation}
where
%
%e3.8 #&#
%
\begin{equation}
\label{eq:Lambda} \Lambda= \Sigma_Y^{-1/2} {
\Sigma_{YX}}\Sigma_X^{-1} {\Sigma_{XY}}
\Sigma_Y^{-1/2} \in\R^{q \times q}.
\end{equation}
\end{theorem}

\begin{pf}
We may assume, with no loss of generality, that $\mu$ is the zero
vector. Since ${\Sigma_{X}}$ and
${\Sigma_{Y}}$ both are positive definite the inverse square-roots,
$\Sigma_X^{-1/2}$ and $\Sigma_Y^{-1/2}$, exist.

By considering the standardized variables $\widetilde{X} =
\Sigma_X^{-1/2} X$ and $\widetilde{Y} = \Sigma_Y^{-1/2} Y$, we may
replace the covariance matrix $\Sigma$ by
\[
\widetilde\Sigma = \lleft( %
\matrix{ I_p & {
\Lambda_{XY}}
\cr
{\Lambda_{XY}}' & I_q
} %
\rright),
\]
where
%
%e3.9 #&#
%
\begin{equation}
\label{eq:LambdaXY} {\Lambda_{XY}}= \Sigma_X^{-1/2} {
\Sigma_{XY}}\Sigma_Y^{-1/2}.
\end{equation}
Once we have made these reductions, it follows that the matrix
$\Lambda$ in (\ref{eq:Lambda}) can be written as $\Lambda= {\Lambda_{XY}}'
{\Lambda_{XY}}$ and that it has norm less than or equal to $1$.
Indeed, by the
\emph{partial Iwasawa decomposition} of $\widetilde\Sigma$, viz.,
the identity,
\[
\widetilde\Sigma= \lleft( %
\matrix{ I_p & 0
\cr
{
\Lambda_{XY}}' & I_q } %
\rright) \lleft( %
\matrix{ I_p & 0
\cr
0 &
I_q - {\Lambda_{XY}}'{\Lambda_{XY}}
} %
\rright) \lleft( %
\matrix{ I_p & {
\Lambda_{XY}}
\cr
0 & I_q } %
\rright),
\]
where the zero matrix of any dimension is denoted by $0$, we see
that the matrix $\widetilde{\Sigma}$ is positive semidefinite
if and only if $I_q - \Lambda$ is positive semidefinite. Hence,
$\Lambda\leq I_q$ in the Loewner ordering and therefore $\| \Lambda\|
\leq1$.

We proceed to calculate the distance covariance $\widetilde{\V}(X,Y) =
\V(\widetilde{X},\widetilde{Y})$. It is well-known that the
characteristic function of $(\widetilde{X},\widetilde{Y})$ is
\[
f_{\widetilde{X},\widetilde{Y}}(s,t) = \exp \biggl[ -\frac{1}2 \left( %
\matrix{
s
\cr
t } %
\right)' \widetilde\Sigma \left( %
\matrix{ s
\cr
t } %
\right) \biggr] = \exp \biggl[ -\frac{1}2
\bigl(|s|_p^2 + |t|_q^2 +
2s' {\Lambda_{XY}}t\bigr) \biggr],
\]
where $s \in\R^p$ and $t \in\R^q$. Therefore,
\[
\bigl | f_{\widetilde{X},\widetilde{Y}}(s,t) - f_{\widetilde{X}}(s) f_{\widetilde{Y}}(t)
\bigr |^2 = \bigl( 1 - \exp\bigl(-s' {\Lambda_{XY}}t
\bigr) \bigr)^2 \exp\bigl(-|s|_p^2 -
|t|_q^2\bigr),
\]
and hence
%
%e3.10 #&#
%
\begin{eqnarray}
\label{eq:cpcqvxy} c_p c_q \V^2(\widetilde{X},
\widetilde{Y}) &= & \int_{\R^{p+q}} \bigl( 1 - \exp
\bigl(-s' {\Lambda_{XY}}t\bigr) \bigr)^2 \exp
\bigl(-|s|_p^2-|t|_q^2\bigr)
\frac{\mathrm{d} s}{|s|_p^{p+1}} \frac{\mathrm{d} t}{|t|_q^{q+1}}
\nonumber
\\[-8pt]
\\[-8pt]
&= & \int_{\R^{p+q}} \bigl( 1 - \exp\bigl(s' {
\Lambda_{XY}}t\bigr) \bigr)^2 \exp \bigl(-|s|_p^2
- |t|_q^2\bigr) \frac{\mathrm{d} s}{|s|_p^{p+1}} \frac{\mathrm{d}
t}{|t|_q^{q+1}},
\nonumber
\end{eqnarray}
where the latter integral is obtained by making the change of variables
$s \mapsto-s$ within the former integral.

By a Taylor series expansion, we obtain
\begin{eqnarray*}
\bigl( 1 - \exp\bigl(s' {\Lambda_{XY}}t\bigr)
\bigr)^2 & = & 1 - 2 \exp\bigl(s' {\Lambda_{XY}}t
\bigr) + \exp\bigl(2s' {\Lambda_{XY}}t\bigr)
\\
& = & \sum_{k=2}^\infty\frac{2^k-2}{k!}
\bigl(s' {\Lambda_{XY}}t\bigr)^k.
\end{eqnarray*}
Substituting this series into (\ref{eq:cpcqvxy}) and interchanging
summation and integration, a procedure which is straightforward to
verify by means of Fubini's theorem, and noting that the odd-order
terms integrate to zero, we obtain
%
%e3.11 #&#
%
\begin{equation}
\label{eq:series} c_p c_q \V^2(\widetilde{X},
\widetilde{Y}) = \sum_{k=1}^\infty
\frac{2^{2k}-2}{(2k)!} \int_{\R^{p+q}} \bigl(s' {
\Lambda_{XY}}t\bigr)^{2k} \exp\bigl(-|s|_p^2-|t|_q^2
\bigr) \frac{\mathrm{d} s}{|s|_p^{p+1}} \frac{\mathrm{d} t}{|t|_q^{q+1}}.
\end{equation}

To calculate, for $k \ge1$, the integral
%
%e3.12 #&#
%
\begin{equation}
\label{eq:k} \int_{\R^{p+q}} \bigl(s' {
\Lambda_{XY}}t\bigr)^{2k} \exp\bigl(-|s|_p^2-|t|_q^2
\bigr) \frac{\mathrm{d} s}{|s|_p^{p+1}} \frac{\mathrm{d} t}{|t|_q^{q+1}},
\end{equation}
we change variables to polar coordinates, putting $s = r_x \theta$ and
$t = r_y \phi$ where $r_x, r_y > 0$, $\theta= (\theta_1, \ldots,
\theta_p)' \in S^{p-1}$, and $\phi= (\phi_1, \ldots, \phi_q)' \in
S^{q-1}$. Then the integral (\ref{eq:k}) separates into a product of
multiple integrals over $(r_x,r_y)$, and over $(\theta,\phi)$,
respectively. The integrals over $r_x$ and $r_y$ are standard gamma
integrals,
%
%e3.13 #&#
%
\begin{equation}
\label{eq:k.factor1} \int_0^\infty\int
_0^\infty r_x^{2k-2}
r_y^{2k-2} \exp\bigl(-r_x^2-r_y^2
\bigr) \dd r_x \dd r_y = \tfrac14  \bigl[\Gamma \bigl(k-
\tfrac{1}{2} \bigr) \bigr]^2 =  \bigl[ \bigl(-\tfrac{1}{2}
 \bigr)_k  \bigr]^2 \pi,
\end{equation}
and the remaining factor is the integral
%
%e3.14 #&#
%
\begin{equation}
\label{eq:k.factor2} \int_{S^{q-1}} \int_{S^{p-1}}
\bigl(\theta' {\Lambda_{XY}}\phi\bigr)^{2k} \dd
\theta \dd\phi,
\end{equation}
where $\dd\theta$ and $\dd\phi$ are unnormalized surface measures on
$S^{p-1}$ and $S^{q-1}$, respectively. By a standard invariance argument,
\[
\int_{S^{p-1}} \bigl(\theta'v\bigr)^{2k}
\dd\theta = |v|_p^{2k} \int_{S^{p-1}}
\theta_1^{2k} \dd\theta,
\]
$v \in\R^p$. Setting $v = {\Lambda_{XY}}\phi$ and applying some well-known
properties of the surface measure $\dd\theta$, we obtain
\begin{eqnarray*}
\int_{S^{p-1}} \bigl(\theta' {\Lambda_{XY}}
\phi\bigr)^{2k} \dd\theta & = & |{\Lambda_{XY}}
\phi|_p^{2k} \int_{S^{p-1}}
\theta_1^{2k} \dd \theta
\\
& = & 2 c_{p-1} \frac{\Gamma(k+\trup{1}{2})\Gamma(\trup{1}{2} p)} {
\Gamma(k+(\trup{1}{2}) p)\Gamma(\trup{1}{2})} \bigl(\phi' \Lambda\phi
\bigr)^k.
\end{eqnarray*}
Therefore, in order to evaluate (\ref{eq:k.factor2}), it remains to
evaluate
\[
J_k(\Lambda) = \int_{S^{q-1}} \bigl(
\phi' \Lambda\phi\bigr)^k \dd\phi.
\]
Since the surface measure is invariant under transformation
$\phi\mapsto K\phi$, $K \in O(q)$, it follows that
$J_k(\Lambda) = J_k(K'\Lambda K)$ for all $K \in O(q)$.
Integrating with respect to the normalized Haar measure on the
orthogonal group, we conclude that
%
%e3.15 #&#
%
\begin{equation}
\label{eq:ilambdaintegral} J_k(\Lambda) = \int_{O(q)}
J_k\bigl(K'\Lambda K\bigr) \dd K = \int
_{S^{q-1}} \int_{O(q)} \bigl(
\phi'K'\Lambda K\phi\bigr)^k \dd K \dd \phi.
\end{equation}
We now use the properties of the zonal polynomials. By
(\ref{eq:property.d}),
\[
\bigl(\phi'K'\Lambda K\phi\bigr)^k =
\bigl(\tr K'\Lambda K\phi\phi'\bigr)^k =
\sum_{|\kappa|=k} C_\kappa\bigl(K'
\Lambda K\phi\phi'\bigr);
\]
therefore, by (\ref{eq:property.e}),
\[
\int_{O(q)} \bigl(\phi'K'\Lambda K
\phi\bigr)^k \dd K = \sum_{|\kappa|=k} \int
_{O(q)} C_\kappa\bigl(K'\Lambda K\phi
\phi'\bigr) \dd K = \sum_{|\kappa|=k}
\frac{C_\kappa(\Lambda) C_\kappa(\phi\phi')}{C_\kappa(I_q)}.
\]
Since $\phi\phi'$ is of rank $1$ then, by property (c),
$C_\kappa(\phi\phi') = 0$ if $\ell(\kappa) > 1$; it now follows, by
(\ref{eq:property.d}) and the fact that $\phi\in S^{q-1}$, that
\[
C_{(k)}\bigl(\phi\phi'\bigr) = \sum
_{|\kappa|=k} C_\kappa\bigl(\phi\phi'\bigr) =
\bigl(\tr\phi\phi'\bigr)^k = \bigl(\phi'\phi
\bigr)^k = |\phi|_q^{2k} = 1.
\]
Therefore,
\[
\int_{O(q)} \bigl(\phi'K'\Lambda K
\phi\bigr)^k \dd K = \frac{C_{(k)}(\Lambda)}{C_{(k)}(I_q)} = \frac{(\trup{1}{2})_k}{((\trup{1}{2}) q)_k}
C_{(k)}(\Lambda),
\]
where the last equality follows by (\ref{eq:Ck}). Substituting
this result at (\ref{eq:ilambdaintegral}), we obtain
\[
J_k(\Lambda) = 2 c_{q-1} \frac{(\trup{1}{2})_k}{((\trup{1}{2}) q)_k}
C_{(k)}(\Lambda).
\]
Collecting together these results, and using the well-known identity
$(2k)! = k! 2^{2k} (\trup{1}{2})_k$, we obtain the representation
(\ref{eq:aidcov.mvn}), as desired.
\end{pf}

We remark that by interchanging the roles of $X$ and $Y$ in Theorem~\ref{th:aidcov}, we would obtain (\ref{eq:aidcov.mvn}) with $\Lambda$
in (\ref{eq:Lambda}) replaced by
\[
\Lambda_0 = \Sigma_X^{-1/2} {
\Sigma_{XY}}\Sigma_Y^{-1} {\Sigma_{YX}}
\Sigma_X^{-1/2} \in\R^{p \times p}.
\]
Since $\Lambda$ and $\Lambda_0$ have the same characteristic polynomial
and hence the same set of nonzero eigenvalues, and noting that
$C_\kappa(\Lambda)$ depends only on the eigenvalues of $\Lambda$, it
follows that $C_{(k)}(\Lambda) = C_{(k)}(\Lambda_0)$. Therefore, the
series representation (\ref{eq:aidcov.mvn}) for $\widetilde{\V}^2(X,Y)$
remains unchanged if the roles of $X$ and $Y$ are interchanged.

The series appearing in Theorem~\ref{th:aidcov} can be expressed
in terms of the generalized hypergeometric functions of matrix argument
(Gross and Richards \cite{GroRic87}, James \cite{Jam64}, Muirhead
\cite{Mui82}). For
this purpose, we introduce the \emph{partitional rising factorial} for
any $\alpha\in\C$ and any partition $\kappa= (k_1,\ldots,k_q)$ as
\[
(\alpha)_\kappa= \prod_{j=1}^q
\bigl(\alpha-(\trup {1} {2}) (j-1) \bigr)_{k_j}.
\]
Let $\alpha_1,\ldots,\alpha_l,\beta_1,\ldots,\beta_m \in\C$ where
$-\beta_i + \frac{1}2(j-1)$ is not a nonnegative integer, for all
$i = 1,\ldots,m$ and $j = 1,\ldots,q$. Then the ${}_lF_m$ generalized
hypergeometric function of matrix argument is defined as\vspace*{-1pt}
\[
{}_lF_m(\alpha_1,\ldots,
\alpha_l;\beta_1,\ldots,\beta_m;S) = \sum
_{k=0}^{\infty} \frac{1}{k!} \sum
_{|\kappa|=k} \frac{(\alpha_1)_\kappa\cdots(\alpha_l)_\kappa} {
(\beta_1)_\kappa\cdots(\beta_m)_\kappa} C_\kappa(S),
\]
where $S$ is a symmetric matrix. A complete analysis of the
convergence properties of this series was derived by Gross and
Richards \cite{GroRic87}, page 804, Theorem~6.3, and we refer the
reader to that
paper for the details.

%co3.2 #&#
%
\begin{corollary} \label{cor:aidcov}
In the setting of Theorem~\ref{th:aidcov}, we have\vspace*{-1pt}
%
%e3.16 #&#
%
\begin{eqnarray}
\label{eq:aidcov} \widetilde{\V}^2(X,Y) &=& 4 \pi\frac{c_{p-1}}{c_p}
\frac{ c_{q-1}}{c_q} \biggl( {}_3F_2 \biggl(\frac{1}2,-
\frac{1}2,-\frac{1}2;\frac{1}2 p,\frac{1}2 q;
\Lambda \biggr)
\nonumber
\\[-8pt]
\\[-8pt]
&&\phantom{4 \pi\frac{c_{p-1}}{c_p}
\frac{ c_{q-1}}{c_q} \biggl(}{} - 2 {}_3F_2 \biggl(\frac{1}2,-
\frac{1}2,-\frac{1}2;\frac{1}2 p,\frac{1}2 q;
\frac{1}4\Lambda \biggr) + 1 \biggr).
\nonumber
\end{eqnarray}
\end{corollary}

\begin{pf}
It is evident that
\[
(\trup{1} {2})_\kappa= %
\cases{ (\trup{1}
{2})_{k_1}, &\quad\mbox{if } $\ell(\kappa) \le1$,
\cr
0, &\quad\mbox{if
} $\ell(\kappa) > 1$.} %
\]
Therefore, we now can write the series in (\ref{eq:aidcov.mvn}),
up to a multiplicative constant, in terms of a generalized
hypergeometric function of matrix argument, in that
\begin{eqnarray*}
&&\sum_{k=1}^\infty \frac{2^{2k}-2}{k! 2^{2k}}
\frac{(\trup{1}{2})_k (-\trup{1}{2})_k (-\trup{1}{2})_k} {
((\trup{1}{2}) p)_k ((\trup{1}{2}) q)_k} C_{(k)}(\Lambda)
\\
&&\quad = \sum_{k=1}^\infty\frac{2^{2k}-2}{k! 2^{2k}}
\sum_{|\kappa|=k} \frac{(\trup{1}{2})_\kappa(-\trup{1}{2})_\kappa(-\trup
{1}{2})_\kappa} {
((\trup{1}{2}) p)_\kappa((\trup{1}{2}) q)_\kappa} C_\kappa (
\Lambda)
\\
&&\quad= \sum_{k=1}^\infty\frac{1}{k!}
\sum_{|\kappa|=k} \frac{(\trup{1}{2})_\kappa(-\trup{1}{2})_\kappa(-\trup
{1}{2})_\kappa} {
((\trup{1}{2}) p)_\kappa((\trup{1}{2}) q)_\kappa} C_\kappa (
\Lambda)
\\
&&\qquad{} - 2 \sum_{k=1}^\infty
\frac{1}{k! 2^{2k}} \sum_{|\kappa|=k} \frac{(\trup{1}{2})_\kappa(-\trup{1}{2})_\kappa(-\trup
{1}{2})_\kappa} {
((\trup{1}{2}) p)_\kappa((\trup{1}{2}) q)_\kappa}
C_\kappa (\Lambda)
\\
&&\quad = \biggl[ {}_3F_2 \biggl(\frac{1}2,-
\frac{1}2,-\frac{1}2;\frac{1}2 p,\frac{1}2 q;
\Lambda \biggr) - 1 \biggr] - 2 \biggl[ {}_3F_2 \biggl(
\frac{1}2,-\frac{1}2,-\frac{1}2;\frac{1}2 p,
\frac{1}2 q;\frac{1}4\Lambda \biggr) - 1 \biggr].
\end{eqnarray*}
Due to property (\ref{eq:homogeneous}) it remains to show that the
zonal polynomial series expansion for the
${}_3F_2  (\frac{1}2,-\frac{1}2,-\frac{1}2;\frac{1}2 p,\frac{1}2 q;
\Lambda )$
generalized hypergeometric function of matrix argument converges
absolutely for all $\Lambda$ with $\Lambda\leq I_q$ in the Loewner
ordering. By (\ref{eq:Ck})
\begin{eqnarray*}
{}_3F_2 \biggl(\frac{1}2,-\frac{1}2,-
\frac{1}2;\frac{1}2 p,\frac{1}2 q; \Lambda \biggr) &
\leq&\sum_{k=0}^\infty\frac{2^{2k}}{k! 2^{2k}}
\frac{(-\trup{1}{2})_k (-\trup{1}{2})_k} {
((\trup{1}{2}) p)_k }
\\
&=& {}_2F_1 \biggl(-\frac{1}2,-
\frac{1}2;\frac{1}2 p;1 \biggr).
\end{eqnarray*}
The latter series converges due to Gauss' theorem for hypergeometric
functions and so we have absolute convergence
at (\ref{eq:aidcov}) for all $\Sigma$ with positive definite marginal
covariance
matrices.
\end{pf}

Consider the case in which $q = 1$ and $p$ is arbitrary. Then
$\Lambda$ is a scalar; say, $\Lambda= \rho^2$ for some
$\rho\in[-1,1]$. Then the ${}_3F_2$ generalized hypergeometric
functions in (\ref{eq:aidcov}) each reduce to a Gaussian
hypergeometric function, denoted by ${}_2F_1$, and
(\ref{eq:aidcov}) becomes
\[
\widetilde{\V}^2(X,Y) = 4 \frac{c_{p-1}}{c_p} \biggl(
{}_2F_1 \biggl(-\frac{1}2,-\frac{1}2;
\frac{1}2 p;\rho^2 \biggr) - 2 {}_2F_1
\biggl(-\frac{1}2,-\frac{1}2;\frac{1}2 p;
\frac{1}4\rho^2 \biggr) + 1 \biggr).
\]
For the case in which $p = q = 1$, we may identify $\rho$ with the
Pearson correlation coefficient and the hypergeometric series can be
expressed in terms of elementary functions. By well-known results
(Andrews, Askey and Roy \cite{AndAskRoy99}, pages 64 and 94),
%
%e3.17 #&#
%
\begin{equation}
\label{eq:hgfidentity} {}_2F_1 \bigl(-\tfrac{1}2,-
\tfrac{1}2;\tfrac{1}2;\rho^2\bigr) = \rho
\sin^{-1} \rho+ \bigl(1-\rho^2\bigr)^{1/2},
\end{equation}
and thus we derive the same result for $p = q = 1$ as in Sz\'ekely
\textit{et al.} \cite{SzeRizBak07}, page 2785.

For cases in which $q = 1$ and $p$ is odd, we
can again obtain explicit expressions for $\widetilde{\V}^2(X,Y)$.
In such cases, the ${}_3F_2$ generalized hypergeometric functions
in (\ref{eq:aidcov}) reduce to Gaussian hypergeometric functions
of the form ${}_2F_1 (-\frac{1}2,-\frac{1}2;k+\frac{1}2;\rho^2)$,
$k \in\N$, and it can be shown that these latter functions are
expressible in closed form in terms of elementary functions and the
$\sin^{-1}(\cdot)$ function. For instance, for $p = 3$, the contiguous
relations for the ${}_2F_1$ functions can be used to show that\vspace*{1pt}
%
%e3.18 #&#
%
\begin{equation}
\label{eq:oddp} {}_2F_1 \biggl(-\frac{1}2,-
\frac{1}2;\frac{3}2;\rho^2\biggr) =
\frac{3(1-\rho^2)^{1/2}}{4} + \frac{(1+2\rho^2)\sin^{-1}\rho}{4\rho}.
\end{equation}
Further, by repeated application of the same contiguous relations,
it can be shown that for $k = 2,3,4,\ldots$\,,
\[
{}_2F_1 \bigl(-\tfrac12,-\tfrac12;k+\tfrac12;
\rho^2\bigr) = \rho^{-2(k-1)} \bigl(1-\rho^2
\bigr)^{1/2} P_{k-1}\bigl(\rho^2\bigr) +
\rho^{-(2k-1)} Q_k\bigl(\rho^2\bigr)
\sin^{-1} \rho,
\]
where $P_k$ and $Q_k$ are polynomials of degree $k$.
Therefore, for $q = 1$ and $p$ odd, the distance
covariance $\widetilde{\V}^2(X,Y)$ can be expressed in closed form
in terms of elementary functions and the $\sin^{-1}(\cdot)$ function.

The appearance of the generalized hypergeometric functions of matrix
argument also yields a useful expression for the affinely invariant
distance variance. In order to state this result, we shall define
for each positive integer $p$ the quantity\vspace*{1pt}
%
%e3.19 #&#
%
\begin{equation}
\label{eq:Ap} A(p) = \frac{\Gamma((\trup{1}{2}) p) \Gamma((\trup{1}{2}) p + 1)} {
[\Gamma ((\trup{1}{2})(p+1) ) ]^2} - 2 {}_2F_1
\biggl(-\frac{1}2,-\frac{1}2;\frac{1}2 p;
\frac{1}4 \biggr) + 1.
\end{equation}

%co3.3 #&#
%
\begin{corollary} \label{cor:aidvar}
In the setting of Theorem~\ref{th:aidcov}, we have\vspace*{1pt}
%
%e3.20 #&#
%
\begin{equation}
\label{eq:aidvar} \widetilde{\V}^2(X,X) = 4 \pi\frac{c_{p-1}^2}{c_{p}^2} A(p).
\end{equation}
\end{corollary}

\begin{pf}
We are in the special case of Theorem~\ref{th:aidcov} for which $X =
Y$, so that $p = q$ and $\Lambda= I_p$. By applying (\ref{eq:Ck}),
we can write the series in (\ref{eq:aidcov.mvn}) as\vspace*{1pt}
\begin{eqnarray*}
&&4 \pi\frac{c_{p-1}^2}{c_{p}^2} \sum_{k=1}^\infty
\frac{2^{2k}-2}{k! 2^{2k}} \frac{(\trup{1}{2})_k (-\trup{1}{2})_k (-\trup{1}{2})_k} {
((\trup{1}{2}) p)_k ((\trup{1}{2}) p)_k} C_{(k)}(I_p)
\\
&&\quad= 4 \pi\frac{c_{p-1}^2}{c_{p}^2} \sum_{k=1}^\infty
\frac{2^{2k}-2}{k! 2^{2k}} \frac{ (-\trup{1}{2})_k (-\trup{1}{2})_k}{((\trup{1}{2}) p)_k}
\\
&&\quad= 4 \pi\frac{c_{p-1}^2}{c_{p}^2} \biggl( \biggl[ {}_2F_1
\biggl(-\frac{1}2,-\frac{1}2;\frac{1}2 p;1 \biggr) - 1
\biggr]
\\
&&\qquad{} - 2 \biggl[ {}_2F_1 \biggl(-\frac{1}2,-
\frac{1}2;\frac{1}2 p;\frac{1}4 \biggr) - 1 \biggr]
\biggr).
\end{eqnarray*}
By Gauss' theorem for hypergeometric functions the series
${}_2F_1(-\frac{1}2,-\frac{1}2;\frac{1}2 p;z)$ also converges for the
special value $z = 1$, and then
\[
{}_2F_1 \biggl(-\frac{1}2,-\frac{1}2;
\frac{1}2 p;1\biggr) = \frac{\Gamma((\trup{1}{2}) p) \Gamma((\trup{1}{2}) p + 1)} {
[\Gamma ((\trup{1}{2})(p+1) ) ]^2},
\]
thereby completing the proof.
\end{pf}

For cases in which $p$ is odd, we can proceed as explained at
(\ref{eq:oddp}) to obtain explicit values for the Gaussian hypergeometric
function remaining in (\ref{eq:aidvar}). This leads in such cases to
explicit expressions for the exact value of $\widetilde{\V}^2(X,X)$.
In particular, if $p = 1$ then it follows from (\ref{eq:cp}) and
(\ref{eq:hgfidentity}) that
\[
\widetilde{\V}^2(X,X) = \frac{4}3 - \frac{4(\sqrt{3}-1)}{\pi};
\]
and for $p = 3$, we deduce from (\ref{eq:cp}) and (\ref{eq:oddp}) that
\[
\widetilde{\V}^2(X,X) = 2 - \frac{4(3\sqrt{3}-4)}{\pi}.
\]

Corollaries~\ref{cor:aidcov} and~\ref{cor:aidvar} enable the
explicit and efficient calculation of the affinely invariant distance
correlation (\ref{eq:aidcor.def}) in the case of subvectors of a
multivariate normal population. In doing so, we use the algorithm of
Koev and Edelman \cite{KoeEde06} to evaluate the generalized
hypergeometric
function of matrix argument, with C and Matlab code being available at
these authors' websites.

%f1 #&#
%
\begin{figure}%[t]

\includegraphics{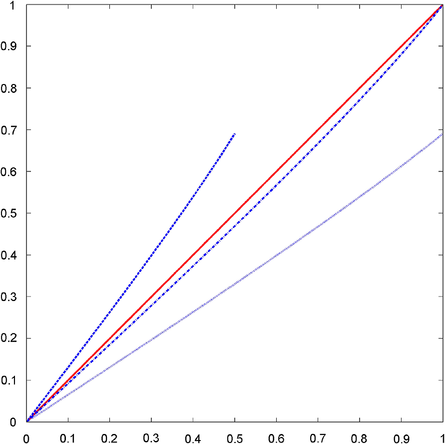}

\caption{The affinely invariant distance correlation for subvectors of a
multivariate normal population, where $p = q = 2$, as a function of
the parameter $r$ in three distinct settings. The solid diagonal
line is the identity function and is provided to serve as a reference
for the three distance correlation functions. See the text for details.}
\label{fig:aidcor.22}
\end{figure}

Figure~\ref{fig:aidcor.22} concerns the case $p = q = 2$ in various
settings, in which the matrix $\Lambda_{22}$ depends on a single
parameter $r$ only. The dotted line shows the affinely invariant
distance correlation when
\[
{\Lambda_{XY}}= %
\pmatrix{ 0 & 0
\cr
0 & r } %
;
\]
this is the case with the weakest dependence considered here. The
dash-dotted line applies when
\[
{\Lambda_{XY}}= %
\pmatrix{ r & 0
\cr
0 & r } %
.
\]
The strongest dependence corresponds to the dashed line, which shows
the affinely invariant distance correlation when
\[
{\Lambda_{XY}}= %
\pmatrix{ r & r
\cr
r & r } %
;
\]
in this case we need to assume that $0 \le r \le\frac{1}2$ in order
to retain positive definiteness.

%f2 #&#
%
\begin{figure}%[t]

\includegraphics{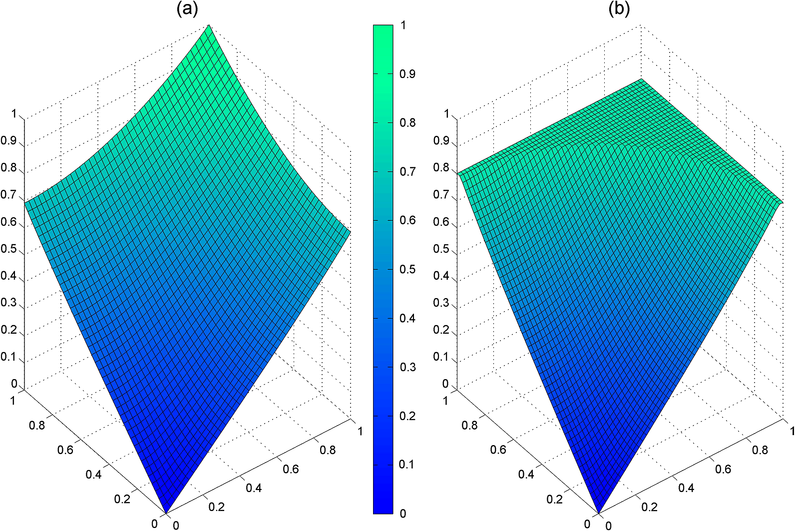}

\caption{The affinely invariant distance correlation between the
$p$- and $q$-dimensional subvectors of a $(p+q)$-dimensional
multivariate normal population, where (a) $p = q = 2$ and
${\Lambda_{XY}}= \diag(r,s)$, and (b)~$p = 2$, $q = 1$ and ${\Lambda
_{XY}}= (r,s)'$.}
\label{fig:aidcor.rs}
\end{figure}

In Figure~\ref{fig:aidcor.rs}, panel (a) shows the affinely invariant
distance correlation when $p = q = 2$ and
\[
{\Lambda_{XY}}= %
\pmatrix{ r & 0
\cr
0 & s } %
,
\]
where $0 \le r, s \le1$. With reference to Figure~\ref
{fig:aidcor.22}, the margins correspond to the dotted line and
the diagonal corresponds to the dash-dotted line.

Panel (b) of Figure~\ref{fig:aidcor.rs} concerns the case in which
$p = 2$, $q = 1$ and ${\Lambda_{XY}}= (r,s)'$, where $r^2 + s^2 \le
1$. Here,
the affinely invariant distance correlation attains an upper limit as
$r^2 + s^2 \uparrow1$, and we have evaluated that limit numerically
as $0.8252$.

%s4 #&#
\section{Limit theorems} \label{sec:limits}

We now study the limiting behavior of the affinely invariant distance
correlation measures for subvectors of multivariate normal populations.

Our first result quantifies the asymptotic decay of the affinely
invariant distance correlation in the case in which the
cross-covariance matrix converges to the zero matrix, in that
\[
\tr(\Lambda) = \| {\Lambda_{XY}}\|_F^2
\longrightarrow0,
\]
where $\| \cdot\|_F$ denotes the Frobenius norm, and the
matrices $\Lambda= {\Lambda_{XY}}' {\Lambda_{XY}}$ and ${\Lambda
_{XY}}$ are defined in
(\ref{eq:Lambda}) and (\ref{eq:LambdaXY}), respectively.

%th4.1 #&#
%
\begin{theorem} \label{th:limit.1}
Suppose that $(X,Y) \sim\mathcal{N}_{p+q}(\mu,\Sigma)$, where
\[
\Sigma= \lleft( %
\matrix{ {\Sigma_{X}}& {
\Sigma_{XY}}
\cr
{\Sigma_{YX}}& {\Sigma_{Y}} }
\rright)
\]
with ${\Sigma_{X}}\in\R^{p \times p}$ and ${\Sigma_{Y}}\in\R^{q
\times q}$ being
positive definite, and suppose that the matrix $\Lambda$ in
(\ref{eq:Lambda})
has positive trace. Then,
%
%e4.1 #&#
%
\begin{equation}
\label{eq:limit.1} \lim_{\tr(\Lambda) \to0} \frac{\widetilde{\mathcal{R}}^2(X,Y)}{\tr(\Lambda)} =
\frac{1}{4 pq \sqrt{A(p)A(q)}},
\end{equation}
where $A(p)$ is defined in (\ref{eq:Ap}).
\end{theorem}

\begin{pf}
We first note that $\widetilde{\V}^2(X,X)$ and $\widetilde{\V}^2(Y,Y)$
do not depend on ${\Sigma_{XY}}$, as can be seen from their explicit
representations in terms of $A(p)$ and $A(q)$ given in (\ref{eq:aidvar}).

In studying the asymptotic behavior of $\widetilde{\V}^2(X,Y)$, we may
interchange the limit and the summation in the series representation
(\ref{eq:aidcov.mvn}). Hence, it suffices to find the limit term-by-term.
Since $C_{(1)}(\Lambda) = \tr(\Lambda) $ then the ratio of the
term for $k = 1$ and $\tr(\Lambda)$ equals
\[
\frac{c_{p-1}}{c_p} \frac{c_{q-1}}{c_q} \frac{\pi}{pq}.
\]
For $k \geq2$, it follows from (\ref{eq:zonal.onepart}) that
$C_{(k)}(\Lambda)$ is a sum of monomials in the eigenvalues of
$\Lambda$, with each monomial being of degree $k$, which is greater
than the degree, viz. $1$, of $\tr(\Lambda)$; therefore,
\[
\lim_{\tr(\Lambda) \to0} \frac{C_{(k)}(\Lambda)}{\tr(\Lambda)} = \lim_{\Lambda\to0}
\frac{C_{(k)}(\Lambda)}{\tr(\Lambda)} = 0.
\]
Collecting these facts together, we obtain (\ref{eq:limit.1}).
\end{pf}

If $p = q = 1$, we are in the situation of Theorem~7(iii) in Sz\'ekely
\textit{et al.} \cite{SzeRizBak07}. Applying the identity (\ref{eq:hgfidentity}), we obtain
\[
{}_2F_1\biggl(-\frac{1}2,-\frac{1}2;
\frac{1}2;\frac{1}4\biggr) %= \frac12 \sin^{-1} \frac12 + \frac{\sqrt{3}}{2}
= \frac{\pi}{12} +
\frac{\sqrt{3}}{2},
\]
and $ (\tr(\Lambda))^{1/2} = |\rho|$. Thus, we obtain
\[
\lim_{\rho\to0} \frac{\widetilde{\mathcal{R}}(X,Y)}{|\rho|} = \frac{1}{2  (1 + (\trup{1}{3}) \pi- \sqrt{3} )^{1/2}},
\]
as shown by Sz\'ekely \textit{et al.} \cite{SzeRizBak07}, page 2785.

In the remainder of this section, we consider situations in which one
or both of the dimensions $p$ and $q$ grow without bound. We will
repeatedly make use of the fact that, with $c_p$ defined as in
(\ref{eq:cp}),
%
%e4.2 #&#
%
\begin{equation}
\label{eq:lem.1} \frac{c_{p-1}}{\sqrt{p} c_p} \longrightarrow\frac{1}{\sqrt{2\pi}}
\end{equation}
as $p \to\infty$, which follows easily from the functional
equation for the gamma function along with Stirling's formula.

%th4.2 #&#
%
\begin{theorem} \label{th:limit.2}
For each positive integer $p$, suppose that
$(X_p,Y_p) \sim\mathcal{N}_{2p}(\mu_p,\Sigma_p)$, where
\[
\Sigma_p = \lleft( %
\matrix{ {\Sigma_{X,p}} & {
\Sigma_{XY,p}}
\cr
{\Sigma_{YX,p}} & {\Sigma_{Y,p}} }
\rright)
\]
with ${\Sigma_{X,p}} \in\R^{p \times p}$ and ${\Sigma_{Y,p}}
\in\R^{p \times p}$ being positive definite and such that
\[
\Lambda_p = \Sigma^{-1/2}_{Y, p} {
\Sigma_{YX,p}} \Sigma_{X, p}^{-1} {\Sigma_{XY,p}}
\Sigma^{-1/2}_{Y, p} \neq0.
\]
Then
%
%e4.3 #&#
%
\begin{equation}
\label{eq:limit.2a} \lim_{p \to\infty} \frac{p}{\tr(\Lambda_p)} \widetilde{
\mathcal{V}}^2(X_p,Y_p) = \frac{1}{2}
\end{equation}
and
%
%e4.4 #&#
%
\begin{equation}
\label{eq:limit.2b} \lim_{p \to\infty} \frac{p}{\tr(\Lambda_p)} \widetilde{
\mathcal{R}}^2(X_p,Y_p) = 1.
\end{equation}
\end{theorem}

In particular, if $\Lambda_p = r^2 I_p$ for some $r \in[0,1]$ then
$\tr(\Lambda_p) = r^2p$, and so (\ref{eq:limit.2a}) and
(\ref{eq:limit.2b}) reduce to
\[
\lim_{p \to\infty} \widetilde{\V}^2(X_p,Y_p)
= \tfrac12 r^2 \quad\mbox{and}\quad \lim_{p \to\infty}
\widetilde{\mathcal{R}}(X_p,Y_p) = r,
\]
respectively. The following corollary concerns the special case in
which $r = 1$; we state it separately for emphasis.

%co4.3 #&#
%
\begin{corollary} \label{cor:limit.2}
For each positive integer $p$, suppose that
$X_p \sim\mathcal{N}_{p}(\mu_p,\Sigma_p)$,
with $\Sigma_p$ being positive definite. Then
%
%e4.5 #&#
%
\begin{equation}
\label{eq:limit.2c} \lim_{p \to\infty} \widetilde{\V}^2(X_p,X_p)
= \tfrac12.
\end{equation}
\end{corollary}

\begin{pf*}{Proof of Theorem~\ref{th:limit.2} and Corollary~\ref{cor:limit.2}} In order to prove (\ref{eq:limit.2a}), we study
the limit for the terms corresponding separately to $k=1$, $k=2$, and
$k \ge3$ in (\ref{eq:aidcov.mvn}).

For $k=1$, on recalling that $C_{(1)}(\Lambda_p)=\tr(\Lambda_p)$,
it follows from (\ref{eq:lem.1}) that the ratio of that term to
$\tr(\Lambda_p) / p$ tends to $1/2$.

For $k=2$, we first deduce from (\ref{eq:property.d}) that
$C_{(2)}(\Lambda_p) \le(\tr\Lambda_p)^2$. Moreover,
$\tr(\Lambda_p) \le p$ because $\Lambda_p \le I_p$ in the Loewner
ordering. Thus, the ratio of the second term in (\ref{eq:aidcov.mvn})
to $\tr(\Lambda_p) / p$ is a constant multiple of
\begin{eqnarray*}
\frac{p}{\tr(\Lambda_p)} \frac{c_{p-1}^2}{c_p^2} \frac{C_{(2)}(\Lambda_p)}{((\trup{1}{2}) p)_2 ((\trup{1}{2}) p)_2} &\le&\frac{c_{p-1}^2}{c_p^2}
\frac{p^2}{((\trup{1}{2}) p)_2 ((\trup{1}{2}) p)_2}
\\
&=& 4 \frac{p}{(p+1)^2} \frac{c_{p-1}^2}{p c_p^2}
\end{eqnarray*}
which, by (\ref{eq:lem.1}), converges to zero as $p \to\infty$.

Finally, suppose that $k \ge3$. Obviously,
$\Lambda_p \le\| \Lambda_p \| I_p$ in the Loewner ordering inequality,
and so it follows from (\ref{eq:zonal.onepart}) that
$C_{(k)}(\Lambda_p) \le\| \Lambda_p \|^k C_{(k)}(I_p)$.
Also, since $\tr(\Lambda_p) \ge\| \Lambda_p \|$ then by again
applying the Loewner ordering inequality and (\ref{eq:Ck}) we
obtain
%
%e4.6 #&#
%
\begin{equation}
\label{eq:limit.inequalities} \frac{C_{(k)}(\Lambda_p)}{\tr(\Lambda_p)} \le \frac{\|\Lambda_p\|^k C_{(k)}(I_p)}{\| \Lambda_p \|} = \|
\Lambda_p \|^{k-1} C_{(k)}(I_p) \le
C_{(k)}(I_p) = \frac{((\trup{1}{2}) p)_k}{(\trup{1}{2})_k}.
\end{equation}
Therefore,
\begin{eqnarray*}
&&4\pi\frac{p}{\tr(\Lambda_p)} \frac{c_{p-1}^2}{c_p^2} \sum_{k=3}^\infty
\frac{2^{2k}-2}{k! 2^{2k}} \frac{(\trup{1}{2})_k (-\trup{1}{2})_k (-\trup{1}{2})_k} {
((\trup{1}{2}) p)_k ((\trup{1}{2}) p)_k} C_{(k)}(\Lambda_p)
\\
&&\quad\le4\pi p \frac{c_{p-1}^2}{c_p^2} \sum_{k=3}^\infty
\frac{2^{2k}-2}{k! 2^{2k}} \frac{(-\trup{1}{2})_k (-\trup{1}{2})_k}{((\trup{1}{2}) p)_k}.
\end{eqnarray*}
By\vspace*{1pt} (\ref{eq:lem.1}), each term $pc_{p-1}^2/(\frac{1}2 p)_k c_p^2$
converges to zero as $p \to\infty$, and this proves both
(\ref{eq:limit.2a}) and its special case, (\ref{eq:limit.2c}).
Then, (\ref{eq:limit.2b}) follows immediately.
\end{pf*}

Finally, we consider the situation in which $q$, the dimension of $Y$,
is fixed while $p$, the dimension of $X$, grows without bound.

%th4.4 #&#
%
\begin{theorem} \label{th:limit.3}
For each positive integer $p$, suppose that $(X_p,Y) \sim
\mathcal{N}_{p+q}(\mu_p,\Sigma_p)$, where
\[
\Sigma_p = \lleft( %
\matrix{ {\Sigma_{X,p}} & {
\Sigma_{XY,p}}
\cr
{\Sigma_{YX,p}} & {\Sigma_{Y}} }
\rright)
\]
with ${\Sigma_{X,p}} \in\R^{p \times p}$ and ${\Sigma_{Y}}\in
\R^{q \times q}$
being positive definite and such that
\[
\Lambda_p = \Sigma^{-1/2}_Y {
\Sigma_{YX,p}} \Sigma_{X, p}^{-1} {\Sigma_{XY,p}}
\Sigma^{-1/2}_Y \neq0.
\]
Then
%
%e4.7 #&#
%
\begin{equation}
\label{eq:limit.3a} \lim_{p \to\infty} \frac{\sqrt{p}}{\tr(\Lambda_p)} \widetilde{
\mathcal{V}}^2(X_p,Y) = \sqrt{\frac{\pi}{2}}
\frac{c_{q-1}}{q c_q }
\end{equation}
and
%
%e4.8 #&#
%
\begin{equation}
\label{eq:limit.3b} \lim_{p \to\infty} \frac{\sqrt{p}}{\tr(\Lambda_p)} \widetilde{
\mathcal{R}}^2(X_p,Y) = \frac{1}{2q \sqrt{A(q)}}.
\end{equation}
\end{theorem}

\begin{pf}
By (\ref{eq:aidcov.mvn}),
\[
\widetilde{\V}^2(X_p,Y) = 4\pi
\frac{c_{p-1}}{c_p} \frac{c_{q-1}}{c_q} \sum_{k=1}^\infty
\frac{2^{2k}-2}{k! 2^{2k}} \frac{(\trup{1}{2})_k (-\trup{1}{2})_k (-\trup{1}{2})_k} {
((\trup{1}{2}) p)_k ((\trup{1}{2}) q)_k} C_{(k)}(\Lambda_p).
\]
We now examine the limiting behavior, as $p \to\infty$, of
the terms in this sum for $k=1$ and, separately, for $k \ge2$.

For $k=1$, the limiting value of the ratio of the corresponding term to
$\tr(\Lambda_p)/\sqrt{p}$ equals
\[
\pi\frac{c_{q-1}}{q c_q} \lim_{p \to\infty} \frac{\sqrt{p}}{\tr(\Lambda_p)}
\frac{c_{p-1}}{p c_p} C_{(1)}(\Lambda_p) = \sqrt{
\frac{\pi}{2}} \frac{c_{q-1}}{q c_q }
\]
by (\ref{eq:lem.1}) and the fact that
$C_{(1)}(\Lambda_p)= \tr(\Lambda_p)$.

For $k \geq2$, the ratio of the sum to $\tr(\Lambda_p)/\sqrt{p}$
equals
\begin{eqnarray*}
&&4 \pi\frac{\sqrt{p}}{\tr(\Lambda_p)} \frac{c_{p-1}}{c_p} \frac
{c_{q-1}}{c_q} \sum
_{k=2}^\infty \frac{2^{2k}-2}{k! 2^{2k}} \frac{(\trup{1}{2})_k (-\trup{1}{2})_k (-\trup{1}{2})_k}{((\trup{1}{2}) p)_k
((\trup{1}{2}) q)_k}
C_{(k)}(\Lambda_p)
\\
&&\quad \leq4 \pi\frac{\sqrt{p}}{\|\Lambda_p\|} \frac{c_{p-1}}{c_p} \frac{c_{q-1}}{c_q} \sum
_{k=2}^\infty\frac{2^{2k}-2}{k! 2^{2k}}
\frac{(-\trup{1}{2})_k (-\trup{1}{2})_k}{((\trup{1}{2}) p)_k} \|\Lambda_p\|^k
\\
&&\quad \leq4 \pi\sqrt{p} \frac{c_{p-1}}{c_p} \frac{c_{q-1}}{c_q} \sum
_{k=2}^\infty\frac{2^{2k}-2}{k! 2^{2k}} \frac{(-\trup{1}{2})_k (-\trup{1}{2})_k}{((\trup{1}{2}) p)_k},
\end{eqnarray*}
where we have used (\ref{eq:limit.inequalities}) to obtain the last
two inequalities. By applying (\ref{eq:lem.1}), we see that the latter
upper bound converges to $0$ as $p \to\infty$, which proves (\ref
{eq:limit.3a}),
and then (\ref{eq:limit.3b}) follows immediately.
\end{pf}

The results in this section have practical implications for
affine distance correlation analysis of large-sample, high-dimensional
Gaussian data. In the setting of Theorem~\ref{th:limit.3},
$\tr(\Lambda_p) \leq q$ is bounded, and so
\[
\lim_{p \to\infty} \widetilde{\mathcal{R}}(X_p,Y) = 0.
\]
As a consequence of Theorem~\ref{th:consistency} on the consistency of
sample measures, it follows that the direct calculation of affine
distance correlation measures for such data will return values which
are virtually zero. In practice, in order to obtain values of the
sample affine distance correlation measures which permit statistical
inference, it will be necessary to calculate $\widehat{\Lambda}_p$,
the maximum likelihood estimator of $\Lambda_p$, and then to rescale
the distance correlation measures with the factor
$\sqrt{p}/\tr(\widehat{\Lambda}_p)$. In the scenario of Theorem~\ref{th:limit.2}, the asymptotic behavior of the affine distance
correlation measures depends on the ratio $p/\tr(\Lambda_p)$; and as
$\tr(\Lambda_p)$ can attain any value in the interval $[0,p]$, a wide
range of asymptotic rates of convergence is conceivable.

In all these settings, the series representation (\ref{eq:aidcov.mvn})
can be used to obtain complete asymptotic expansions in powers of
$p^{-1}$ or $q^{-1}$, of the affine distance covariance or correlation
measures, as $p$ or $q$ tend to infinity.

%s5 #&#
\section{Time series of wind vectors at the Stateline wind energy center}
\label{sec:Stateline}

R{\'e}millard \cite{Rem09} proposed the use of the distance
correlation to
explore nonlinear dependencies in time series data. Zhou
\cite{Zho12}
pursued this approach recently and defined the auto distance covariance
function and the auto distance correlation function, along with
natural sample versions, for a strongly stationary vector-valued time
series, say $(X_j)_{j = - \infty}^\infty$.

It is straightforward to
extend these notions to the affinely invariant distance correlation.
Thus, for an integer $k$, we refer to
%
%e5.1 #&#
%
\begin{equation}
\label{eq:acf} \widetilde{\mathcal{R}}_X(k) = \frac{\widetilde{\V}(X_j,X_{j+k})}{\widetilde{\V}(X_j,X_j)}
\end{equation}
as the \emph{affinely invariant auto distance correlation} at the lag
$k$. Similarly, given jointly strongly stationary, vector-valued time
series $(X_j)_{j= -\infty}^{\infty}$ and $(Y_j)_{j=
-\infty}^{\infty}$, we refer to
%
%e5.2 #&#
%
\begin{equation}
\label{eq:ccf} \widetilde{\mathcal{R}}_{X,Y}(k) = \frac{\widetilde{\V}(X_j,Y_{j+k})} {
\sqrt{\widetilde{\V}(X_j,X_j) \widetilde{\V}(Y_j,Y_j)}}
\end{equation}
as the \emph{affinely invariant cross distance correlation} at the
lag $k$. The corresponding sample versions can be defined in the natural
way, as in the case of the non-affine distance correlation (Zhou \cite{Zho12}).

We illustrate these concepts on time series data of wind observations
at and near the Stateline wind energy center in the Pacific Northwest
of the United States. Specifically, we consider time series of bivariate
wind vectors at the meteorological towers at Vansycle, right at the
Stateline wind farm at the border of the states of Washington and
Oregon, and at Goodnoe Hills, 146 km west of Vansycle along the
Columbia River Gorge. Further information can be found in the paper
by Gneiting \textit{et al.} \cite{Gneetal06}, who developed a regime-switching space-time
(RST) technique for 2-hour-ahead forecasts of hourly average wind
speed\vadjust{\goodbreak} at the Stateline wind energy center, which was then the largest
wind farm globally. For our purposes, we follow Hering and Genton \cite
{HerGen10}
in studying the time series at the original 10-minute resolution, and
we restrict our analysis to the longest continuous record, the 75-day
interval from August 14, 2002 to October 28, 2002.

%f3 #&#
%
\begin{figure}%[t]

\includegraphics{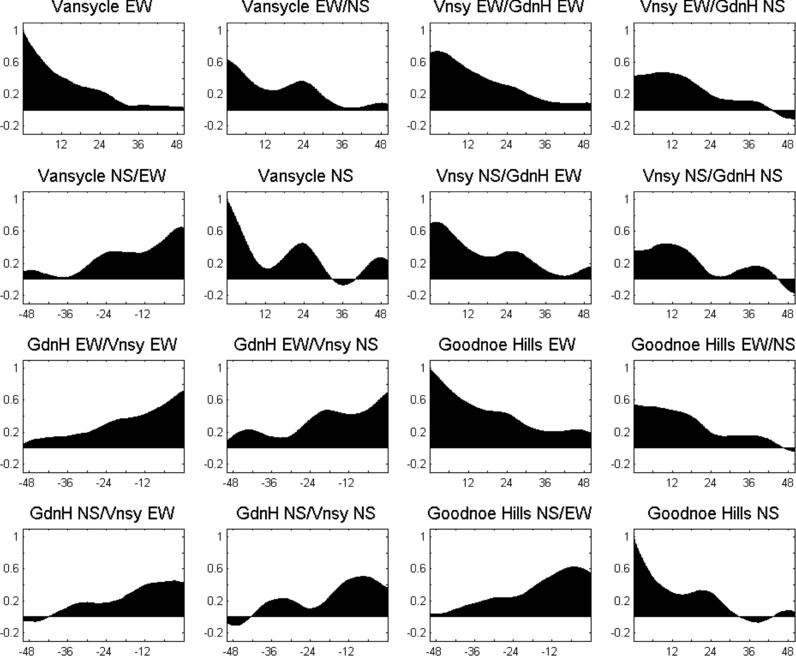}

%{\includegraphics[height=360pt, width=480pt]{Pearson4x4.png}}
\caption{Sample auto and cross Pearson correlation functions for the
univariate time series $V^{\mathrm{EW}}_j$, $V^{\mathrm{NS}}_j$,
$G^{\mathrm{EW}}_j$, and $G^{\mathrm{NS}}_j$, respectively. Positive
lags indicate
observations at the westerly site (Goodnoe Hills) leading those at
the easterly site (Vansycle), or observations of the north--south
component leading those of the east--west component, in units of
hours.}
\label{fig:4x4pcorplot}
\end{figure}

%f4 #&#
%
\begin{figure}%[t]

\includegraphics{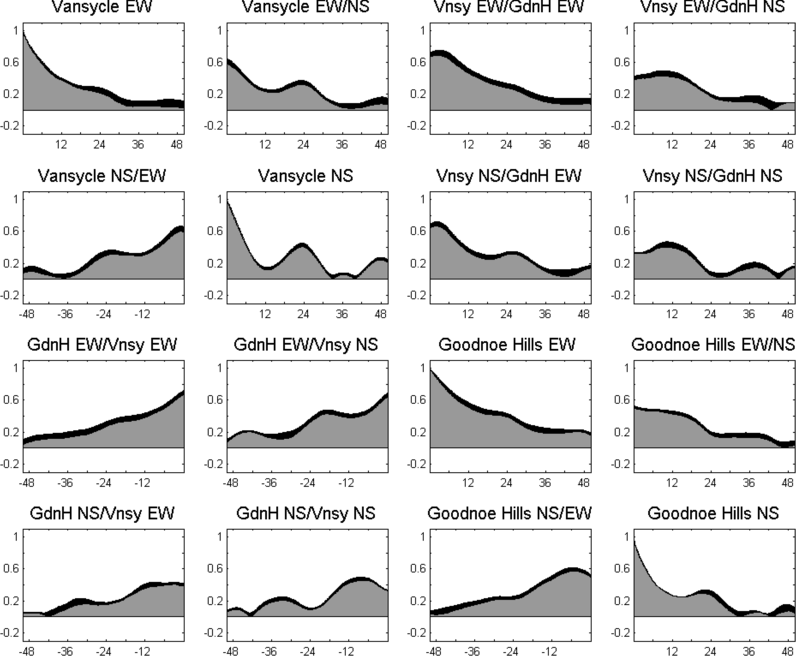}

%{\includegraphics[height=360pt, width=480pt]{comparison4x4.png}}
\caption{Sample auto and cross distance correlation functions for the
univariate time series $V^{\mathrm{EW}}_j$, $V^{\mathrm{NS}}_j$,
$G^{\mathrm{EW}}_j$, and $G^{\mathrm{NS}}_j$, respectively. For
comparison, we also
display, in grey, the values that arise when the sample Pearson
correlations in Figure~\protect\ref{fig:4x4pcorplot} are converted to
distance correlation under the assumption of Gaussianity; these
values generally are smaller than the original ones. Positive lags
indicate observations at Goodnoe Hills leading those at Vansycle, or
observations of the north--south component leading those of the
east--west component, in units of hours.}
\label{fig:4x4dcorplot}
\end{figure}

Thus, we consider time series of bivariate wind vectors over $10\,800$
consecutive $10$-minute intervals. We\vspace*{2pt}
write $V^{\mathrm{NS}}_j$ and
$V^{\mathrm{EW}}_j$ to denote the north--south and the east--west component,
respectively, of the wind vector at Vansycle at time $j$, with positive
values corresponding to northerly and easterly winds. Similarly, we
write $G^{\mathrm{NS}}_j$ and $G^{\mathrm{EW}}_j$ for the north--south
and the
east--west component, respectively, of the wind vector at Goodnoe Hills
at time $j$.

Figure~\ref{fig:4x4pcorplot} shows the classical (Pearson) sample auto
and cross correlation functions for the four univariate time series.
The auto correlation functions generally decay with the temporal, but
do so non-monotonously, due to the presence of a diurnal component.
The cross correlation functions between the wind vector components at
Vansycle and Goodnoe Hills show remarkable asymmetries and peak at
positive lags, due to the prevailing westerly and southwesterly wind
(Gneiting \textit{et al.} \cite{Gneetal06}). In another interesting feature, the cross
correlations between the north--south and east--west components at lag
zero are strongly positive, documenting the dominance of southwesterly
winds.

Figure~\ref{fig:4x4dcorplot} shows the sample auto and cross distance
correlation functions for the four time series; as these variables are
univariate, there is no distinction between the standard and the
affinely invariant version of the distance correlation. The patterns
seen resemble those in the case of the Pearson correlation. For
comparison, we also display values of the distance correlation based
on the sample Pearson correlations shown in Figure~\ref
{fig:4x4pcorplot}, and converted to distance correlation under the
assumption of bivariate Gaussianity, using the results of Sz\'ekely
\textit{et al.} \cite{SzeRizBak07}, page 2785, and Section~\ref
{sec:aidc.mvn}; in every single
case, these values are smaller than the original ones.\vadjust{\goodbreak}

%f5 #&#
%
\begin{figure}%[t]

\includegraphics{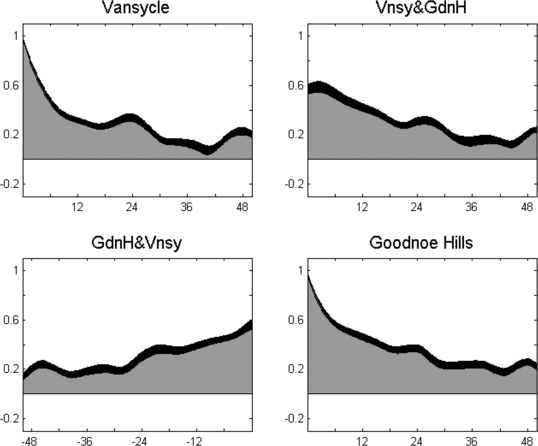}

%{\includegraphics[height=200pt]{2x2vergleichsplot.png}}
\caption{Sample auto and cross affinely invariant distance correlation
functions for the bivariate time series $(V^{\mathrm{EW}}_j,V^{\mathrm
{NS}}_j)'$ and $(G^{\mathrm{EW}}_j,G^{\mathrm{NS}}_j)'$ at Vansycle
and Goodnoe
Hills. For comparison, we also display, in grey, the values that
are generated when the Pearson correlation in Figure~\protect\ref
{fig:4x4pcorplot} is converted to the affinely invariant
distance correlation under the assumption of Gaussianity; these
converted values generally are smaller than the original ones.
Positive lags indicate observations at Goodnoe Hills leading those
at Vansycle, in units of hours.}
\label{fig:2x2corplot}
\end{figure}

Having considered the univariate time series setting, it is natural
and complementary to look at the wind vector time series $(V^{\mathrm
{EW}}_j,V^{\mathrm{NS}}_j)$ at Vansycle and $(G^{\mathrm
{EW}}_j,G^{\mathrm{NS}}_j)$ at
Goodnoe Hills\vspace*{2pt}
 from a genuinely multivariate perspective. To this end,
Figure~\ref{fig:2x2corplot} shows the sample affinely invariant auto
and cross distance correlation functions for the bivariate wind vector
series at the two sites. Again, a diurnal component is visible, and
there is a remarkable asymmetry in the cross-correlation functions,
which peak at lags of about two to three hours.

In light of our analytical results in Section~\ref{sec:aidc.mvn}, we
can compute the affinely invariant distance correlation between
subvectors of a multivariate normally distributed random vector. In
particular, we can compute the affinely invariant auto and
cross distance correlation between bivariate subvectors
of a 4-variate Gaussian process with Pearson auto and cross correlations
as shown in Figure~\ref{fig:4x4pcorplot}. In Figure~\ref{fig:2x2corplot},
values of the affinely invariant distance correlation that have been
derived from Pearson correlations in these ways are shown in grey; the
differences from those values that are computed directly from the data
are substantial, with the converted values being smaller, possibly
suggesting that assumptions of Gaussianity may not be appropriate for
this particular data set.\looseness=1

We wish to emphasize that our study is purely exploratory: it is
provided for illustrative purposes and to serve as a basic example.
In future work, the approach hinted at here may have the potential to be
developed into parametric or nonparametric bootstrap tests for Gaussianity.
For this purpose recall that, in the Gaussian setting, the affinely invariant
distance correlation is a function of the canonical correlation coefficients,
that is, $\widetilde{\mathcal{R}}=g(\lambda_1,\ldots,\lambda_r)$.
For a parametric bootstrap test, one could generate $B$ replicates of
$g(\lambda_1^{\star},\ldots,\lambda_r^{\star})$, leading to a pointwise
($1-\alpha$)-confidence band. The test would now reject Gaussianity
if the sample affinely invariant distance correlation function does not
lie within this
band. For the nonparametric bootstrap test, one could obtain ensembles
$\widetilde{\mathcal{R}}_n^\star$ by resampling methods, again
defining a
pointwise ($1-\alpha$)-confidence band and checking if
$g(\widehat{\lambda}_1,\ldots,\widehat{\lambda}_r)$ is located
within this band.

Following the pioneering work of Zhou \cite{Zho12},
the distance correlation may
indeed find a wealth of applications in exploratory and inferential
problems for time series data.

%s6 #&#
\section{Discussion} \label{sec:discussion}

In this paper, we have studied an affinely invariant version of the
distance correlation measure introduced by Sz\'ekely \textit{et al.}
\cite{SzeRizBak07}
and Sz{\'e}kely and Rizzo \cite{SzeRiz09} in both population and
sample settings
(see Sz{\'e}kely and Rizzo \cite{SzeRiz12} for further aspects of
the role of
invariance in properties of distance correlation measures).
The affinely invariant distance correlation shares the desirable
properties of the standard version of the distance correlation and
equals the latter in the univariate case. In the multivariate case,
the affinely invariant distance correlation remains unchanged under
invertible affine transformations, unlike the standard version, which
is preserved under orthogonal transformations only. Furthermore, the
affinely invariant distance correlation admits an exact and readily
computable expression in the case of subvectors from a multivariate
normal population. We have shown elsewhere that the standard version
allows for a series expansion too, but this does not appear to be a
series that generally can be made simple, and further research will be
necessary to make it accessible to efficient numerical computation.
Related asymptotic results can be found in Gretton \textit{et al.}
\cite{Greetal12} and
Sz{\'e}kely and Rizzo~\cite{SzeRiz13}.\looseness=1

Competing measures of dependence also have featured prominently recently
(Reshef \textit{et al.} \cite{autokey15}, Speed \cite{Spe11}). However,
those measures are
restricted to univariate settings, and claims of superior performance
in exploratory data analysis have been disputed (Gorfine, Heller and Heller
\cite{autokey4}, Simon and Tibshirani \cite{autokey17}).
We therefore share much of
Newton's \cite{New09} enthusiasm about the use of the distance
correlation as
a measure of dependence and association. A potential drawback for large
data sets is the computational cost required to compute the sample
distance covariance, and the development of computationally efficient
algorithms or subsampling techniques for doing this is highly
desirable.\looseness=1

% zodis "Acknowledgments" paliekamas pagal autoriu
\section*{Acknowledgements}

The research of Johannes Dueck, Dominic Edelmann and Tilmann Gneiting
has been supported by the \textit{Deutsche Forschungsgemeinschaft}
(German Research Foundation) within the programme ``Spatio/Temporal
Graphical Models and Applications in Image Analysis,'' grant GRK 1653.

%suskaldyti doi

% imsref loaded by audrone.aklyte, 2014-01-14 11:21:20
% imsref loaded by audrone.aklyte, 2014-01-14 11:40:05
% imsref loaded by audrone.aklyte, 2014-01-14 11:42:13
%

\printhistory

\end{document}